\documentclass[11pt,a4paper]{article}
\usepackage{amsmath}
\usepackage{amsfonts}
\usepackage{mathrsfs}
\usepackage{amsthm,latexsym}
\def\C{\mathbb{C}}
\def\br{\mathbb{R}}
\def\L{\mathcal{L}}
\def\H{\mathbb{H}}
\def\b{\mathscr{B}}
\def\l{\lambda}
\def\r{\mathscr{R}}
\theoremstyle{definition}\newtheorem{thm}{Theorem}
\theoremstyle{definition}
\theoremstyle{definition}\newtheorem{col}{Corollary}
\theoremstyle{definition}\newtheorem{lem}{lemma}
\theoremstyle{remark}
\newcommand{\beq}{\begin{equation}}
\newcommand{\eeq}{\end{equation}}
\newcommand{\beqs}{\begin{eqnarray*}}
\newcommand{\eeqs}{\end{eqnarray*}}
\newcommand{\beqn}{\begin{eqnarray}}
\newcommand{\eeqn}{\end{eqnarray}}
\newcommand{\beqa}{\begin{array}}
\newcommand{\eeqa}{\end{array}}

\begin{document}
\title{\bf Wiener measure for Heisenberg group\footnote{Supported by National Natural Science Foundation of
China under Grant \#10990012 and 50-class General Financial Grant from the China Postdoctoral Science Foundation \#2011M501317.}}
\author{\small Heping Liu$^1$ and Yingzhan Wang$^2$\\[-0.1cm]
\small $^1$ LMAM, School of Mathematical Sciences, Peking University, Beijing 100871, P. R. China\\[-0.1cm]
\small $^2$ School of Sciences, South China University of Technology, Guangzhou 510641, P. R. China\\[-0.1cm]
\small $^2$ Corresponding Author\\[0.15cm]
\small E-mail: $^1$hpliu@pku.edu.cn, $^2$wyzde@pku.edu.cn  \\[0.15cm]}
\date{}
\maketitle \noindent {\bf Abstract.}In this paper, we build Wiener
measure for the path space on the Heisenberg group by using of the
heat kernel corresponding to the sub-Laplacian and give the
definition of the Wiener integral. Then we give the Feynman-Kac
formula. \vskip 0.2cm \noindent {\bf Keywords:} Heisenberg group,
 C-C distance, sub-Laplacian
operator, Wiener measure, Feynman-Kac formula.

\section{Introduction}Wiener measure was first built by using of the heat kernel corresponding
to Laplacian $\triangle$ by Wiener in 1920s. It gives a perfect
explanation to some physics problem. (see [\ref{kac}]). The
connection between Wiener measure and Schr\"{o}dinger equation has
been known for a long time. Kac  in 1949 first gave the explicit
solution expression of the following equation \beq\label{1}\left\{
    \beqa{ll}
      \frac{\partial}{\partial t}u(t,\xi)=\frac{1}{2}\frac{\partial^{2}}{\partial{\xi^{2}}}u(t,\xi)-V(\xi)u(t,\xi)\\
     u(0,\xi)=f(\xi)
    \eeqa
  \right.
\eeq through Wiener integral. And the upper Wiener integral solution
is the famous  Feynman-Kac formula.
Over the recent years, there are still much work on Wiener measure. See \cite{Ryu},\cite{hsu},\cite{sad}.
In this paper, we study the Wiener space on Heisenberg group.

We suppose $\H$  Heisenberg group with underlining manifold $\br^{2n+1}\cong \C^{n}\times \br$.
The multiplication is given by
$$(z,u)(z',u')=(z+z', u+u'+2\Im(z\bar{z'})) .$$
There is a homogeneous norm
$|\cdot|$ satisfying the subadditivity inequality.
$$|\xi|=(|z|^{4}+u^{2})^{\frac{1}{4}},~  \textrm{for}~ \xi=(z,u)\in \H .$$
Obviously it satisfies: $|\xi\eta|\leq |\xi|+|\eta|,~\textrm{and}~
|\xi|=|\xi^{-1}|~\textrm{for}~ \xi, \eta\in \H$. Thus we can define the  distance $N$ on $\H$.
$$N(x,y)=|xy^{-1}|, x,y \in \H.$$
 Let $C_{o}[0,1]$ denote the set of
continuous functions $x(t)$ in the unit interval [0,1] valued in
H with $x(0)=o$, where $o$ is the unit element of
$\H$. Then $C_{o}[0,1]$ is a Banach space with the norm
$\|x\|=\sup_{0\leq t\leq 1}|x(t)|$. And  $C_{o}[0,1]$ is called the Wiener space in Heisenberg group.

If
$X_{i},Y_{i},U (1\leq i \leq n)$ denote the left invariant vector
fields on H whose values at o are given by
$\frac{\partial}{\partial x_{i}},\frac{\partial}{\partial
y_{i}},\frac{\partial}{\partial u} (1\leq i \leq n)$ respectively, then
$$X_{i}=\frac{\partial}{\partial
x_{i}}+2y_{i}\frac{\partial}{\partial
u},Y_{i}=\frac{\partial}{\partial
y_{i}}-2x_{i}\frac{\partial}{\partial u},U=\frac{\partial}{\partial
u}.$$
Let $\L$ be the Kohn-laplacian operator sometimes called sub-Laplacian operator.$$\L=-\sum_{i=1}^{n}(X_{i}^{2}+Y_{i}^{2}) .$$
For $\L$, the heat kernel's explicit expression was first  given by A.  Hulanicki and  B.  Gaveau.
$$p_{t}(x,y,u)=(2\pi)^{-1}(4\pi)^{-n}\int_{\br}(\frac{|\l|}{\sinh|\l|})^{n}\exp\{-\frac{|\l||z|^{2}}{4}\coth|\l|s-i\l\cdot u\}d\l .$$
(see \cite{Hulanicki}, \cite{Gaveau}.)
Then we can  use it to build the wiener measure in Heisenberg group.

 This paper is organized in three sections. In the second section, we build the Wiener measure for the  Wiener space in Heisenberg group. Then we define the Wiener integral.  In the last section, we use it to give the Feynman-Kac formula corresponding to  $\L$. Throughout the whole paper, constant c may not be the same at every appearance.

\section{ Wiener measure and Wiener integral}
As described in the introduction,  Wiener space is the following set,
$$C_{o}[0,1]=\{x: [0,1]\rightarrow \H, \text{continuous}, x(0)=o .\} $$
\noindent It is a Banach space with the norm $\| \cdot \|$. The  corresponding Borel field is denoted by $\b$. To build the Wiener measure, first we will  build a measure on one of its generating algebras. Next we   extend the measure to the whole field.

For E a Borel set in $\mathbb{R}^{2n+1}$, we define the cylinder set  in $C_{o}[0,1]$,
$$I=\left\{x\in C_{o}[0,1]; (x_{t_{1}},x_{t_{2}},...,x_{t_{m}})\in E
\right\} $$ where $0< t_{1}<t_{2}<t_{3}<...<t_{m}\leq 1$. And  denote the
set of all  the cylinder sets  by $\mathscr{R}$.  Obviously $\r$ is a subfield of $\b$.
In fact, we have the following lemma.
\begin{lem}\label{lem1}
$\sigma(\mathscr{R})=\mathscr{B}$
\end{lem}
\begin{proof} Since $\r$ is a subfield of $\b$, so
$$\sigma(\mathscr{R})\subseteq \mathscr{B}.$$
And we can see obviously
\begin{equation*}
\{x\in C_{o}[0,1]\big{|}~ \Arrowvert x\Arrowvert \leq
1\}=\bigcap_{m=1}^{\infty}\{x\in C_{o}[0,1]\big{|}~
|x(\frac{k}{2^{m}})|\leq 1, k=1,2,...,2^{m}\} .
\end{equation*}
So the closed unit ball is in $\sigma(\mathscr{R})$. Then
$$\sigma(\mathscr{R})\supseteq \mathscr{B} .$$  And we  proved the
$\sigma$-field generated by $\mathscr{R}$ is just $\mathscr{B}$.
\end{proof}

For $I,E,p_{t}(\cdot)$  described above, Denote
\begin{equation}\label{2}W(I)=\int_{E}\prod_{j=1}^{m}p_{t_{j}-t_{j-1}}(u_{j-1}^{-1}u_{j})\prod_{j=1}^{m}du_{j}
.
\end{equation}\label{14}
 Then we give our main theorem of this paper.
 \begin{thm} \label{the1}$W$ defined above is a measure
on $\mathscr{R}$. Furthermore  $W$ is $\sigma$-additive on
$\mathscr{R}$. Then  it can be extended to a measure on the whole
$\sigma$-field $\mathscr{B}$.
 \end{thm}
\begin{proof}
Denote by
$$I_{1}=\left\{x\in C_{o}[0,1]:
(x_{t_{1_{1}}},x_{t_{1_{2}}},...,x_{t_{1_{p}}})\in B_{1}\times
B_{2}\times...\times B_{p} \right\} ,
$$
$$I_{2}=\left\{x\in
C_{o}[0,1]: (x_{t_{2_{1}}},x_{t_{2_{2}}},...,x_{t_{2_{q}}})\in
C_{1}\times C_{2}\times...\times C_{q} \right\} ,$$
 where
$B_{i},C_{j},i=1,2,...p~;~j=1,2,...q$, are Borel sets in
$\mathbb{R}^{2n+1}$. To prove $W$ is a measure on $\mathscr{R}$, we
need only prove  when $I_{1}\cap I_{2}=\emptyset$, we have
$$W(I_{1}\bigcup I_{2})=W(I_{1})+W(I_{2}).$$ We make time take the
same values in the upper two sets. To do this, we need only  add some $t_{k}$
and $B_{k}( \text{or}~ C_{k})$, where $B_{k}(\text{or}~
 C_{k})=\mathbb{R}^{2n
 +1}$. So that  we get:
$$I_{1}'=\left\{x\in C_{o}[0,1]:
(x_{t_{1}},x_{t_{2}},...,x_{t_{m}})\in B_{1}\times
B_{2}\times...\times B_{m} \right\},
$$
$$I_{2}'=\left\{x\in C_{o}[0,1]:
(x_{t_{1}},x_{t_{2}},...,x_{t_{m}})\in C_{1}\times
C_{2}\times...\times C_{m} \right\},
$$
Obviously $I_{1}=I_{1}',I_{2}=I_{2}'$. Since $I_{1}\cap
I_{2}=\emptyset$, there must be some $j$ such that $B_{j}\cap
C_{j}=\emptyset$. Then by expression (\ref{14}),
$$W(I_{1}\bigcup I_{2})=W(I_{1}'\bigcup
I_{2}')=W(I_{1}')+W(I_{2}').$$
 So we need only prove this kind of
equality $W(I)=W(I')$, where
\begin{align*} &I=\left\{x\in C_{o}[0,1]:
(x_{t_{1}},x_{t_{2}})\in B_{1}\times B_{2} , t_{1}<t_{2}\right\},\\
&I'=\left\{x\in C_{o}[0,1]:  (x_{t_{1}},x_{t},x_{t_{2}})\in
B_{1}\times\mathbb{R}^{2n+1} \times B_{2}, t_{1}<t<t_{2}\right\}.
\end{align*}
The other cases can be proved analogously. \\
Through the definition of W, and using the properties of the heat kernel, we can get
\begin{align*}
W(I')&=\int_{B_{1}}\int_{\mathbb{R}^{2n+1}}\int_{
B_{2}}p_{t_{1}}(u_{1})p_{t-t_{1}}(u_{1}^{-1}u_{2})p_{t_{2}-t}(u_{2}^{-1}u_{3})du_{1}du_{2}du_{3}\\
&=\int_{B_{1}}\int_{
B_{2}}p_{t_{1}}(u_{1})p_{t_{2}-t_{1}}(u_{1}^{-1}u_{2})du_{1}du_{2}\\
&=W(I).
\end{align*}

  Now we finished the proof that $W$ is a measure on
$\mathscr{R}$. If we can prove $W$ is $\sigma$-additive on
$\mathscr{R}$, then it has an unique extension to $\mathscr{B}$.
We need only prove the following claim :

If $\left\{I_{m}\right\} (m=1,2,...)$ is a
sequence of cylinder sets in $\mathscr{R}$, $I_{j}\supset I_{j+1}$,
and $\cap_{m=1}^{\infty}I_{m}=\emptyset$, then $\lim_{m\rightarrow
\infty}W(I_{m})=0$. \\Let
\begin{eqnarray*}
I_{m}&=&(t^{(m)}_{1},t^{(m)}_{2},...,t^{(m)}_{s_{m}};E_{m})\\
&\equiv&\left\{x\in C_{o}(0,1);
(x^{(m)}_{t_{1}},x^{(m)}_{t_{2}},...,x^{(m)}_{t_{s_{m}}})\in
E_{m}\subset R^{s_{m}(2n+1)}\right\}.
\end{eqnarray*}

For $\forall~ \varepsilon >0$, we choose a closed set $G_{m}\subset
E_{m}$ such that $$W(I_{m}\setminus
K_{m})<\frac{\varepsilon}{2^{m+1}},$$ where
$$K_{m}=(t^{(m)}_{1},t^{(m)}_{2},...,t^{(m)}_{s_{m}};G_{m}).$$ Let
$L_{m}=\bigcap_{j=1}^{m}K_{j}\in \mathscr{R}$. \\
Then $$L_{m}\subset K_{m}\subset I_{m}.$$ Hence
$$W(I_{m})=W(I_{m}\setminus L_{m})+W(L_{m}).$$
Scince $$ I_{m}\setminus L_{m}=I_{m}\setminus
\bigcap_{j=1}^{m}K_{j}=\bigcup_{j=1}^{m}( I_{m}\setminus
K_{j})\subset\bigcup_{j=1}^{m} (I_{j}\setminus K_{j} ),$$ Then
$$W(I_{m}\setminus L_{m})\leq
\sum_{j=1}^{m}\frac{\varepsilon}{2^{m+1}}\leq\frac{\varepsilon}{2}.$$

So we need prove there exists $m_{0}$ such that
$W(L_{m})<\frac{\varepsilon}{2}$ for any $m>m_{0}$. Here we give the following
lemma that we leave its proof behind:
\begin{lem}\label{lem2}
Let $a>0$ and $0<r<\frac{1}{2}$. Denote  $$
H^{r}(a)=\{x\in C_{0}[0,1]| ~|x(t_{2})^{-1}x(t_{1})|\leq
\frac{2a}{1-2^{-r}}|t_{1}-t_{2}|^{r},\forall t_{1},t_{2}\in
[0,1]\}.$$  Then for $I\in \mathscr{R}$, if  $I\subset
H^{r}(a)^{c}$, , we have  $$\lim_{a\rightarrow \infty}W(I)=0.$$
\end{lem}
Now we suppose  Lemma \ref{lem2} is right. If we can prove when $a$ is large enough,
$L_{m}\subset H^{r}(a)^{c}$, then for m large enough,  $W(L_{m})<\frac
{\varepsilon}{2}$ is obious. So the only thing we should prove is there exists
$m_{0}$ such that $M_{m}=L_{m}\bigcap H^{r}(a)=\emptyset$ for any
$m>m_{0}$. Note $\left\{M_{m}\right\}_{m\geq 1}$ is a decreasing
sequence, and $\bigcap_{m=1}^{\infty}M_{m}=\emptyset$. If $M_{m}\neq
\emptyset$ for any $m$, we choose $x_{m}\in M_{m}$. Then
$\left\{x_{m}\right\}_{m\geq 1}$ is equi-continuous and uniformly
bounded. Therefore by Ascoli-Arzelas theorem,
 we can choose
a convergent sequence still denoted by $\left\{x_{m}\right\}_{m\geq
1}$. Suppose $\lim_{m\rightarrow\infty}x_{m}=x_{0} $. Then $x_{0}\in
H^{r}(a)$.  Note that $M_{m}$ is a compact set, and $x_{m}\in M_{m},
m\geq m_{0}$. So $x_{0}\in M_{m},\forall m\geq 1$, which is a
contrary to $\bigcap_{m=1}^{\infty}M_{m}=\emptyset$.  And theorem \ref{the1} is proved. So what left is to prove
Lemma \ref{lem2}.
\end{proof}
\begin{proof}[Proof of Lemma 2]
Because $N(g,h)=|g^{-1}h|$ is a
distance on $\H$, it is subadditive . Then, for $x\in C_{o}[0,1]$, if there exist
$a>0, r>0$ such
that$$|x(\frac{k}{2^{m}})x(\frac{k-1}{2^{m}})^{-1}|\leq
a(\frac{1}{2^{m}})^{r}, \forall~ k=1,2,...2^{m}, \forall m>0,$$ Then
$$|x(t_{1})x(t_{2})^{-1}|\leq \frac{2a}{1-2^{-r}}|t_{1}-t_{2}|^{r},$$ where
$t_{1},t_{2}\in [0,1],~\text{and each of them can be written as}~
\frac{k}{2^{m}}, k  ~\text{is odd}$.

In fact, for $t_{1}\leq t_{2}$, we can choose $t=\frac{q}{2^{p}}$
$q$ is odd, and $p$ the smallest such that $t_{1}\leq t\leq t_{2}$.
If $t\neq t_{1}$, then we can write
$$t-t_{1}=\frac{1}{2^{m_{1}}}+\frac{1}{2^{m_{2}}}+...+\frac{1}{2^{m_{j}}}, m_{1}<m_{2}<...<m_{j}$$
and if $t\neq t_{2}$,
$$t_{2}-t=\frac{1}{2^{n_{1}}}+\frac{1}{2^{n_{2}}}+...+\frac{1}{2^{n_{k}}}, n_{1}<n_{2}<...<n_{k}.$$
Consider the intervals, $$[t_{1},t_{1}+\frac{1}{2^{m_{j}}}],
[t_{1}+\frac{1}{2^{m_{j}}},t_{1}+\frac{1}{2^{m_{j}}}+\frac{1}{2^{m_{j-1}}}],...[t-\frac{1}{2^{m_{1}}},t]$$
and $$[t,t+\frac{1}{2^{n_{1}}}],[t+\frac{1}{2^{n_{1}}},t+\frac{1}{2^{n_{1}}}+\frac{1}{2^{n_{2}}}],
...[t_{2}-\frac{1}{2^{n_{k}}},t_{2}].$$
Let $l=\min\{m_{1},n_{1}\}, s=\max\{m_{j},n_{k}\}$, then we can see
$$|x(t_{1})^{-1}x(t_{2})|\leq
2a\sum_{k=l}^{s}(\frac{1}{2^{k}})^{r}\leq
\frac{2a}{1-2^{-r}}|t_{1}-t_{2}|^{r}.$$ Let
$$I^{r}_{a,k,m}=\left\{x\in
C_{0}(0,1)|~|x(\frac{k}{2^{m}})x(\frac{k-1}{2^{m}})^{-1}|>
a(\frac{1}{2^{m}})^{r})\right\}, k=1,2,...,2^{m}.$$ Then we have
$$H^{r}(a)\supset
\bigcap_{m=0}^{\infty}\bigcap_{k=1}^{2^{m}}(I^{r}_{a,k,m})^{c}.$$ So
$$I\subset \big{(}H^{r}(a)\big{)}^{c}\subset \bigcup_{m=0}^{\infty}\bigcup_{k=1}^{2^{m}}I^{r}_{a,k,m}.$$
It follows that
\begin{eqnarray}\label{3}W(I)\leq\sum_{m=0}^{\infty}\sum_{k=1}^{2^{m}}W(I^{r}_{a,k,m}).
\end{eqnarray}
Denote $$E=\left\{(u_{1},u_{2})\in
\H :~|u_{2}u_{1}^{-1}|>a(\frac{1}{2^m})^{r}\right\},$$
 then
\begin{eqnarray}\label{4}
W(I^{r}_{a,k,m})&=&\int_{E}p_{\frac{k-1}{2^{m}}}(u_{1})p_{\frac{1}{2^{m}}}(u_{2}u_{1}^{-1})du_{1}du_{2}\nonumber \\
&=&\int_{|u_{3}|>a(\frac{1}{2^m})^{r}}\int_{\H}p_{\frac{k-1}{2^{m}}}(u_{1})p_{\frac{1}{2^{m}}}(u_{3})du_{1}du_{3}\nonumber \\
&=&\int_{|u_{3}|>a(\frac{1}{2^m})^{r}}p_{\frac{1}{2^{m}}}(u_{3})du_{3}.
\end{eqnarray}
Now we need an estimate of the heat kernel.  Denote  $d(x,y)$  the Carnot-Carath\'{e}odory distance associated to
$\L$  and $B(x,r)=\{y\in M:
d(x,y)<r)\}$  the corresponding Balls
 .  We refer the reader to \cite{Varopoulos} for those definitions and the following important estimate.
\begin{equation*}
p_{t}(x^{-1}y)\leq
\frac{C_{1}}{\mu(B(x,t^{\frac{1}{2}}))}\exp(-\frac{C_{1}d(x,y))^{2}}{t})),
\forall x,y\in \H.
\end{equation*}
From \cite{adam},  we can  get $d(x,y)\geq N(x,y)$. Then there exists a constant $M$, such that
\begin{equation}\label{5}
p_{t}(\xi)\leq cMt^{-n-1}\exp(-\frac{M^{-1}|\xi|^{2}}{t}), \forall
\xi\in N .
\end{equation}
So
\begin{eqnarray}
\eqref{4}&\leq& cM\int_{|u_{3}|>a(\frac{1}{2^{m}})^{r}}\exp(-M^{-1}2^{m}|u_{3}|^{2})2^{m(n+1)}du_{3}\nonumber\\
&\leq&c\int_{\rho>a(\frac{1}{2^{m}})^{r}}\exp(-M^{-1}2^{m}\rho^{2})2^{m(n+1)}\rho^{2n+1}d\rho\nonumber\\
&\leq&c2^{\frac{m}{2}}\int_{\rho>a(\frac{1}{2^{m}})^{r}}\exp(-M^{-1}2^{m-1}\rho^{2})d\rho\nonumber\\
&\leq&c2^{\frac{m}{2}}\int_{\rho>a(\frac{1}{2^{m}})^{r}}\exp(-M^{-1}2^{m-1}\rho^{2})\frac{\rho}{a(\frac{1}{2^{m}})^{r}}d\rho\nonumber\\
&=&\frac{c}{a}2^{m(r-\frac{1}{2})}e^{-(2M)^{-1}2^{m(1-2r)}a^{2}}\nonumber .
\end{eqnarray} Then by \eqref{3},
\begin{eqnarray*}
W(I)&\leq&\sum_{m=0}^{\infty}\sum_{k=1}^{2^{m}}\frac{c}{a}2^{m(r-\frac{1}{2})}e^{-(2M)^{-1}2^{m(1-2r)}a^{2}}\\
&\leq&\sum_{m=0}^{\infty}\frac{c}{a}2^{m(r+\frac{1}{2})}e^{-(4M)^{-1}m(1-2r)a^{2}}\\
&=&\frac{c}{a}(1-2^{r+\frac{1}{2}}e^{-(4M)^{-1}a^{2}(1-2r)})^{-1}
\end{eqnarray*} Then $$\lim_{a\rightarrow \infty}W(I)=0.$$  The
proof of the lemma \ref{lem2} is finished.
\end{proof}

So $W$ on $\mathscr{R}$ has an unique extension to $\mathscr{B}$
still denoted by $W$. We call $W$ the Wiener measure on
$C_{o}[0,1]$. The correspondent integral is called Wiener integral.
For $f\in C_{o}[0,1]$, we define its  Wiener integral
by
$$E^{W}[f]=\int_{C_{o}[0,1]}f(x)W(dx).$$
$W$ can also be seen as the Wiener measure on the whole function
space $C[0,1]$, through  $\widehat{W}(A):= W(A\bigcap C_{0}[0,1])$,
for any Borel set $A\in C[0,1]$. For  $\xi\in \H$, we
 denote the translation by $(T_{\xi}x)(t):=
\xi x(t),  \forall~ x\in C[0,1]$. Then we can define Wiener measure
$W_{\xi}$ on $C[0,1]$,$W_{\xi}(A):=W(T_{\xi}^{-1}A)$, for any Borel
set $A\in C[0,1]$. The corresponding Wiener integral can be defined similarly.

\section{Feynman-Kac formula}
In this section, we will use Wiener integral to give an explicit solution expression to the following type equation.
\begin{equation}\label{9}\left\{
    \begin{array}{ll}
      (\partial_{t}-\Delta_{z}-4|z|^{2}\partial_{t}^{2}+4\sum_{j=1}^{n}(x_{j}\frac{\partial}{\partial y_{j}}-y_{j}\frac{\partial}{\partial x_{j}}))u(t,\xi)=-V(\xi)u(t,\xi)\\
     u(0,\xi)=f(\xi)
    \end{array}
  \right..
\end{equation}
In the upper equation, let
$$X_{i}=\frac{\partial}{\partial
x_{i}}+2y_{i}\frac{\partial}{\partial
u},Y_{i}=\frac{\partial}{\partial
y_{i}}-2x_{i}\frac{\partial}{\partial u} ,$$
\noindent which are just the left invariant vector fields of  Heisenberg group.
Then equation \eqref{9}, turns into
\begin{equation}\label{10}
\left\{
    \begin{array}{ll}
      (\partial_{t}+\L)u(t,\xi)=-V(\xi)u(t,\xi)\\
     u(0,\xi)=f(\xi)
    \end{array}
  \right..
\end{equation}
Here we follow the classical method.

For almost every  $a\in \H$, we denote $\delta_{a}(d\xi)$ the
dot measure on $\H$. That is, for any measurable function, it
holds that
$$\int_{\H}f(\xi)\delta_{a}(d\xi)=f(a) .$$
Note
$$E^{W}[\delta_{x(t)}(B)]=\int_{B}p_{t}(u)du, ~~~~\text{for any Borel set } ~B ~\text{in} ~\H .$$
So $$\frac{dE^{W}[\delta_{x(t)}(\cdot)]}{d\xi}=p_{t}(\xi).$$ We denote
by $\delta(\xi^{-1}a)$ the density function of $\delta_{a}(d\xi)$.
Denote  $\delta_{t,\xi}(x)=\delta(\xi^{-1}x(t))$ .
Let
$E^{W}[\delta_{t,\xi}(x)]=
\frac{d}{d\xi}E^{W}[\delta_{x(t)}(\cdot)] .$
Then $\delta_{t,\xi}(x)$ has its meaning.

 Now we will give the next two lemmas.
\begin{lem}\label{lem3}
Let  $G(x)$ be a Wiener-integrable function on $C[0,1]$,
 then $E^{W}[G(x)\delta_{x_{(t)}}(d\xi)], (0<t\leq1)$ is a finite measure and absolutely
continuous to the Haar measure on $\H$. Furthermore for any
measurable function on $\H$, we have
\begin{eqnarray}\label{6}
\int_{\H}f(\xi)E^W[G(x)\delta_{x(t)}(d\xi)]=E^{W}[f(x(t))G(x)]
\end{eqnarray}
\end{lem}
\begin{proof}Let $B$ be a Borel set
in $\H$, and $\delta_{x(t)}(B)$  a character function of
cylinder $I=\left\{x| x(t)\in B\right\}$. Obviously
$$E^W[|G(x)\delta_{x(t)}(B)|]\leq E^W[|G(x)|].$$
And if $\mu(B)=0$, then
from the definition of $W$ , we know $W(I)=0$. It follows that
\beqs
E^W[G(x)\delta_{x(t)}(B)]&=&\int G(x)\delta_{x(t)}(B)W(dx)\\
&=&\int_{I}G(x)W(dx)=0
\eeqs
 So
$E^W[G(x)\delta_{x(t)}(d\xi)]$ is absolutely continuous to the Haar
measure on $\H$.

For the proof of formula \eqref{6}, we can see obviously it holds for
simple function. For  common  functions, we can always make a
sequence of simple function convergent to it. Then we can get the
desired result.
\end{proof}
Denote
$E^W[G(x)\delta_{t,\xi}(x)]=\frac{d}{d\xi}E^W[G(x)\delta_{x(t)}(\cdot)].$
And we  have the following lemma.
\begin{lem}\label{lem4}
Let $0<s<t\leq 1, G(x)$  a Wiener-integrable function and $G(x)$
only depends on the value of $x$ on [0,s].Then
\begin{eqnarray}\label{7}
E^W[G(x)\delta_{t,\xi}(x)]=\int_{\H}E^W[\delta_{t-s,\eta^{-1}\xi}(x)]E^W[G(x)\delta_{s,\eta}(x)]d\eta
\end{eqnarray}
\end{lem}
\begin{proof}
By definition, we only need prove for any Borel set $B$ in
$\H$, it holds that
\begin{equation}\label{8}
E^W[G(x)\delta_{x(t)}(B)]=\int_{B}\int_{\H}E^W[\delta_{t-s,\eta^{-1}\xi}(x)]E^W[G(x)\delta_{s,\eta}(x)]d\eta
d\xi.
\end{equation}
Note first that
\begin{equation*}
\frac{d}{d\xi}E^{W}[\delta_{(x(s)^{-1}x(t))}(\cdot)]=E^{W}[\delta_{t-s,\xi}(x)].
\end{equation*}
In fact for any Borel set
$B$,\begin{eqnarray*}E^{W}[\delta_{(x(s)^{-1}x(t))}(B)]&=&\int_{(u_{1}^{-1}u_{2})\in
B}p_{s}(u_{1})p_{t-s}(u_{1}^{-1}u_{2})du_{1}du_{2}\\
&=&\int_{\H}\int_{B}p_{s}(u_{1})p_{t-s}(u_{3})du_{1}du_{3}\\
&=&\int_{B}p_{t-s}(u_{3})du_{3}\\
&=&\int_{B}E^{W}[\delta_{t-s,\xi}(x)]d\xi
\end{eqnarray*}
Then using lemma \ref{lem3},
\begin{eqnarray*}
(right~of~\eqref{8})&=&\int_{\H}\int_{B}E^W[\delta_{t-s,\eta^{-1}\xi}(x)]E^W[G(x)\delta_{s,\eta}(x)]
d\xi d\eta\nonumber\\
&=&\int_{\H} E^W[\chi_{B}(\eta
x(s)^{-1}x(t))]E^W[G(x)\delta_{s,\eta}(x)]
d\eta\nonumber\\
&=&E^W[G(x)\delta_{x(t)}(B)]
\nonumber\\
&=&(left ~of~ \eqref{8}) .
\end{eqnarray*}
\end{proof}
Now we are ready to give the main theorem in this section, Feynman-Kac formula.
\begin{thm}\label{thm2} Let $V(\xi)$  be a  lower bounded integrable
function on the Heisenberg group $\H$, and $f(\xi)$ a bounded
measuable function on $\H$. Then
$$u(t,\xi)=E^{W_{\xi}}[f\big{(}x(t)\big{)}e^{-\int^{t}_{0}V(x(s))ds}]$$ is the
solution of the differential  equation \eqref{9} or \eqref{10}.
\end{thm}
To begin the proof of the
theorem, we should first suppose the following conclusions hold.
\begin{thm}\label{thm3}
Let
$$u(t,\xi)=E^{W}[\delta_{t,\xi}e^{-\int_{0}^{t}V(x(s))ds}].$$
Then $u(t,\xi)$ is the solution of the following  equation
\begin{equation}\label{11}
\left\{
  \begin{array}{ll}
    (\partial_{t}+\L)u(t,\xi)=-V(\xi)u(t,\xi) \\
    u(0,\xi)=\delta_{\xi} \\
   \lim_{\xi\rightarrow\infty}u(t,\xi)=0
  \end{array}
\right..
\end{equation}
\noindent where $\delta_{\xi}=1$, if~$\xi=o$; and 0, otherwise .
\end{thm}
And through this theorem, we can get the following corollary,
\begin{col} In upper partial differential equation, if
$V(\cdot)$ is replaced by $V(\eta~\cdot~)$, Then
$$p(t,\eta,\xi)=u(t,\eta^{-1}\xi)=E^{W}[\delta_{t,\eta^{-1}\xi}(x)e^{-\int_{0}^{t}V(\eta
 x(s))ds}]$$
 is the solution of the following  equation
\begin{equation}\label{12}
\left\{
  \begin{array}{ll}
   (\partial_{t}+\L)p(t,\eta,\xi)=-V(\xi)p(t,\eta,\xi)\\
    p(0,\eta,\xi)=\delta_{\eta^{-1}\xi}\\
   \lim_{\xi\rightarrow\infty}p(t,\eta,\xi)=0
  \end{array}
\right..
\end{equation}
\end{col}
\begin{lem}\label{lem5}
$p(t, \xi, \eta)=p(t,\eta, \xi).$
\end{lem}
\begin{proof}[Proof of the Theorem \ref{thm2}]
From the definition of $W_{\xi}$, for any function $F$ on $C[0,1]$,
the following equality holds,
$$E^{W_{\xi}}[F(x)]=E^{W}[F(\xi x)].$$
at the meaning that one side of integral on the equal sign exists. \\
So
\begin{eqnarray*}
u(t,\xi)&=&E^{W}[f(\xi x(t))e^{-\int_{0}^{t}V(\xi x(s))ds}]\\
&=&\int_{\H}f(\eta)p(t,\xi,\eta)d\eta \\
&=&\int_{\H}f(\eta)p(t,\eta,\xi)d\eta.
\end{eqnarray*}
Obviously $u(t,\xi)$ satisfies equation \eqref{9}.
\end{proof}
\begin{proof}[Proof of Theorem \ref{thm3}]
Since$$e^{-\int_{0}^{t}V(x(s))ds}=1-\int_{0}^{t}V(x(\tau))e^{-\int_{0}^{\tau}V(x(s))ds}d\tau.$$
so
$$
u(t,\xi)=E^{W}[\delta_{t,\xi}(x)]-\int_{0}^{t}E^{W}[V(x(\tau))e^{-\int_{0}^{\tau}V(x(s))ds}\delta_{t,\xi}(x)]d\tau.
$$
 by \eqref{6}, \eqref{7},we get
\begin{eqnarray*}
& &E^{W}[V(x(\tau))e^{-\int_{0}^{\tau}V(x(s))ds}\delta_{t,\xi}(x)]\\
&=&\int_{\H}E^W[\delta_{t-\tau,\eta^{-1}\xi}(x)]E^W[V(x(\tau))e^{-\int_{0}^{\tau}V(x(s))ds}\delta_{\tau,\eta}(x)]d\eta\\
&=&\int_{\H}p_{t-\tau}(\eta^{-1}\xi)E^W[V(x(\tau))e^{-\int_{0}^{\tau}V(x(s))ds}\delta_{\tau,\eta}(x)]d\eta\\
&=&E^{W}[\int_{\H}p_{t-\tau}(\eta^{-1}\xi)V(x(\tau))e^{-\int_{0}^{\tau}V(x(s))ds}\delta_{\tau,\eta}(x)d\eta]\\
&=&E^{W}[p_{t-\tau}((x(\tau)^{-1}\xi)V(x(\tau))e^{-\int_{0}^{\tau}V(x(s))ds}]\\
&=&\int_{\H}p_{t-\tau}(\eta^{-1}\xi)V(\eta)E^W[e^{-\int_{0}^{\tau}V(x(s))ds}\delta_{\tau,\eta}(x)]d\eta\\
&=&\int_{\H}V(\eta)u(\tau,\eta)p_{t-\tau}(\eta^{-1}\xi)d\eta .
\end{eqnarray*}
So
$$u(t,\xi)=p_{t}(\xi)-\int_{0}^{t}\int_{\H}V(\eta)u(\tau,\eta)p_{t-\tau}(\eta^{-1}\xi)d\eta
d\tau .$$ And \eqref{11} follows.
\end{proof}
\begin{proof}[Proof of lemma \ref{lem5}]First we will see
\begin{equation}
p(t,\eta,\xi)=E^{W}[\delta_{t,\eta^{-1}\xi}(x)e^{-\int_{0}^{t}V(\xi
x(t)^{-1}x(t-\tau)d\tau}].
\end{equation}
In fact for any Borel measurable function $f$,\begin{eqnarray*}
&&\int_{\H}f(\eta)p(t,\eta,\xi)d\eta\\&=&\int_{\H}f(\eta)E^{W}[\delta_{t,\eta^{-1}\xi}(x)e^{-\int_{0}^{t}V(\eta
 x(s))ds}]d\eta,\\
&=&\int_{\H}f(\xi\eta^{-1})E^{W}[\delta_{t,\eta}(x)e^{-\int_{0}^{t}V(\xi\eta^{-1}
 x(s))ds}]d\eta\\
&=&E^{W}[\int_{\H}f(\xi\eta^{-1})\delta_{t,\eta}(x)e^{-\int_{0}^{t}V(\xi\eta^{-1}
 x(s))ds}d\eta]\\
&=&E^{W}[f(\xi x(t)^{-1})e^{-\int_{0}^{t}V(\xi x(t)^{-1}
 x(s))ds}]\\
&=&\int_{\H}f(\eta)E^{W}[\delta_{t,\eta^{-1}\xi}(x)e^{-\int_{0}^{t}V(\xi
x(t)^{-1}
 x(s))ds}]d\eta\\
&=&\int_{\H}f(\eta)E^{W}[\delta_{t,\eta^{-1}\xi}(x)e^{-\int_{0}^{t}V(\xi
x(t)^{-1}
 x(t-\tau))d\tau}]d\eta .
\end{eqnarray*}
 Next we only need prove
\begin{equation}\label{13}
E^{W}[\delta_{t,\eta^{-1}\xi}(x)e^{-\int_{0}^{t}V(\xi
x(t)^{-1}x(t-\tau)d\tau}]=E^{W}[\delta_{t,\xi^{-1}\eta}(x)e^{-\int_{0}^{t}V(\xi
x(\tau))d\tau}].
\end{equation}
Define map $T: C_{0}[0,t]\rightarrow C_{0}[0,t],
~Tx(s)=x(t)^{-1}x(t-s)$. Then the left of \eqref{13} is equal to
$E^{W}[\delta_{t,\xi^{-1}\eta}(Tx)e^{-\int_{0}^{t}V(\xi\cdot
(Tx)(\tau))d\tau}]$. So if we can prove Wiener measure is invariant
under transform $T$, formula \eqref{13} will holds. \\
\indent In $C_{0}[0,t]$, we need only consider the cylinder sets of
the following forms$$I=\left\{x\in C_{0}[0,t];
(x(t_{1}),x(t_{2}),...,x(t_{n})~)\in
(B_{1},B_{2},...,B_{n})\right\},$$ where
$0<t_{1}<t_{2}<...<t_{n}=t$, and $B_{1},B_{2},...,B_{n}$ are the
Borel sets in $\br^{2n+1}$.\\ Then
\begin{eqnarray*} TI&=&\left\{
(y(t-t_{1}),y(t-t_{2}),...,y(t)^{-1}~)\in
(x(t)^{-1}B_{1},x(t)^{-1}B_{2},...,B_{n})\right\}\\
&=&\left\{ (y(t-t_{n-1}),y(t-t_{n-2}),...,y(t)~)\in
(y(t)B_{n-1},y(t)B_{n-2},...,B_{n}^{-1})\right\}.
\end{eqnarray*}
We need only prove $W(I)=W(TI)$.
In fact, by definition
\begin{eqnarray*}
W(TI)&=&\int_{B_{n}^{-1}}\int_{u_{n}B_{1}}...\int_{u_{n}B_{n-1}}p_{t-t_{n-1}}(u_{1})p_{t_{n-1}-t_{n-2}}(u_{1}^{-1}u_{2})...p_{t_{1}}(u_{n-1}^{-1}u_{n})du_{1}du_{2}...du_{n}\\
&=&\int_{B_{n}^{-1}}\int_{B_{1}}...\int_{B_{n-1}}p_{t-t_{n-1}}(u_{n}v_{1})p_{t_{n-1}-t_{n-2}}(v_{1}^{-1}v_{2})...p_{t_{1}}(u_{n-1}^{-1})dv_{1}dv_{2}...du_{n}\\
&=&\int_{B_{n}}\int_{B_{1}}...\int_{B_{n-1}}p_{t-t_{n-1}}(v_{n}^{-1}v_{1})p_{t_{n-1}-t_{n-2}}(v_{1}^{-1}v_{2})...p_{t_{1}}(v_{n-1}^{-1})dv_{1}dv_{2}...dv_{n}\\
&=&\int_{B_{n}}\int_{B_{n-1}}...\int_{B_{1}}p_{t-t_{n-1}}(v_{n-1}^{-1}v_{n})p_{t_{n-1}-t_{n-2}}(v_{n-2}^{-1}v_{n-1})...p_{t_{1}}(v_{1})dv_{1}dv_{2}...dv_{n}\\
&=&W(I).
\end{eqnarray*}
Now the proof is finished.
\end{proof}

\end{document}